\documentclass[a4paper,12pt]{article}

\usepackage[latin1]{inputenc} 
\usepackage[T1]{fontenc}  
\usepackage{geometry}        
\usepackage[francais, english]{babel}  
\usepackage{amsmath}
\usepackage{amsfonts}
\usepackage{amssymb}
\usepackage{fancyhdr}
\usepackage{url}
\usepackage{mathrsfs}

\usepackage[all]{xy}

\newtheorem{df}{Definition}[section]
\newtheorem{prop}[df]{Proposition}
\newtheorem{thm}[df]{Theorem}

\newtheorem{conj}{Conjecture}

\newcommand{\prf}{\textit{Proof}}
\newtheorem{rmk}[df]{Remark}

\newcommand{\mom}{\mathcal{PM}_{\Omega}}
\newcommand{\tmom}{\widetilde{\mathcal{PM}}_{\Omega}}

\newcommand{\dbar}{\overline{\partial}}
\newcommand{\ddbar}{\partial\overline{\partial}}

\newcommand{\C}{\mathbb{C}}

\newcommand{\XD}{X\backslash D}

\newcommand{\vol}{\operatorname{vol}}

\newcommand{\vareps}{\varepsilon}

\newcommand{\vl}{\operatorname{Vol}}

\newcommand{\scal}{\mathbf{s}}

\newcommand{\cqfd}{ \hfill $\square$ }

\title{\textbf{Metrics of Poincaré type with constant scalar curvature: a topological constraint}}  
\author{\textsc{Hugues AUVRAY}}%
\date{}

\begin{document}

\makeatletter
\renewcommand%
   {\section}%
   {%
   \@startsection{section}%
      {1}%
      {0mm}%
      {\baselineskip}%
      {0.5\baselineskip}%
      {\sc\large\centering}%
   }%
\makeatother

\makeatletter
\renewcommand%
   {\subsubsection}%
   {%
   \@startsection{subsubsection}%
      {1}%
      {0mm}%
      {1.25\baselineskip}%
      {0.25\baselineskip}%
      {\sf\normalsize}%
   }%
\makeatother

\maketitle

\renewcommand{\abstractname}{Abstract}
 \begin{abstract}
Let $D=\sum_{j=1}^N D_j$ a divisor with simple normal crossings in a Kähler manifold $(X,\omega_0)$ of complex dimension $m\geq2$. 
The purpose of this short note is to show that the existence of a Poincaré type metric $\varpi$ with constant scalar curvature in $\mathcal{PM}_{[\omega_0]}$ on $\XD$ implies for all $j$ the inequality $\overline{\scal}<\overline{\scal_{D_j}}$. 
Here $\overline{\scal}$ denotes the scalar curvature of $\varpi$, whereas $\overline{\scal_{D_j}}$ denotes the mean scalar curvature associated to $\mathcal{PM}_{[\omega_0|_{D_j}]}$ or $[\omega_0|_{D_j}]$, depending on that $D_j$ intersects other components or not. 
We also explain how those results were already conjectured by G. Székelyhidi when $D$ is reduced to one component.
  \end{abstract}

 \begin{center}
   \rule{3cm}{0.25pt}
 \end{center}

~

\section{Introduction}

As canonical geometric objects, Kähler metrics with constant scalar curvature on compact manifolds have intensively been studied over the past few years. 
The question of their uniqueness has been settled by works of S.K. Donaldson \cite{don1}, T. Mabuchi \cite{mab1}, X.X. Chen \cite{chen1} and by Chen and G. Tian \cite{ch-ti}. 

On the other hand, about the existence of Kähler metrics with constant scalar curvature (on compact manifolds), fewer is known. 
Nevertheless, in the projective case, Yau suggested \cite{yau} a conjecture relating the existence of such metrics among a fixed polarization class with algebro-geometric properties of the polarized manifold. 
This conjecture has been reformulated by Tian \cite{tia} and Donaldson \cite{don2}, and can be stated as:
 \begin{conj}[Yau, Tian, Donaldson] \label{conj}
  A compact polarized manifold $(M,L)$ is K-stable if, and only if, there exists a Kähler metric with constant scalar curvature in the class $c_1(L)$.
 \end{conj}
This is still an open question, at least in the "only if" direction. 
Indeed, the implication "existence of metric with constant scalar curvature implies K-stability" has been proved by J. Stoppa and Mabuchi, see \cite{sto, mab2}, and the references within.
Concerning the converse, in has been proved so far for toric surfaces by Donaldson \cite{don3}, but remains a very delicate problem in the general case.

In this respect, the aim of this note is to provide necessary conditions to the existence of a metric with constant scalar curvature among a class of Kähler metrics with cusp singularities along a divisor, which we call metrics of Poincaré type. 
Considering metrics of Poincaré type takes its interest in that they can appear as limits when working with sequence of smooth metrics. 
Studying this singular case may also enlighten the smooth case. 
For example, Donaldson \cite{don4} suggests the study of Kähler metrics with conical singularities as a method toward the understanding of the existence of smooth Kähler-Einstein metrics in the Fano case; indeed conical singularities vanish when letting their angle go to $2\pi$, but we can recall that they also tend to cusp singularities, when the angle goes to 0.

Let us recall briefly the terminology of Poincaré type Kähler metrics, in the sense of \cite[part 1]{auv}. 
We consider a divisor $D$ with simple normal crossings in a compact Kähler manifold $(X,\omega_0,J)$ of dimension $m\geq2$ and write its decomposition into smooth irreducible components as $D=\sum_{j=1}^N D_j$. 
By \textit{simple normal crossings} we mean that every irreducible component is smooth, and that around a point where exactly $k$ components, $D_1,\dots,D_k$ say, intersect, one has an open set $U$ of holomorphic coordinates $(z_1,\dots,z_k,z_{k+1},\dots,z_m)$ such that $D_{\ell}\cap U=\{z_{\ell}=0\}$ for all $\ell=1,\dots,k$.
 
Let us endow each line bundle $[D_j]$ with a smooth hermitian metric $|\cdot|_j$, and denote by $\sigma_j\in\mathcal{O}([D_j])$ a holomorphic section such that $D_j=\{\sigma_j=0\}$, $j=1,\dots,N$. 
Up to multiplying $|\cdot|_j$ by a positive constant or a smooth positive function for those $j$, we can assume that $|\sigma_j|_j^2\leq e^{-1}$ so that $\rho_j:=-\log(|\sigma_j|^2_j)\geq 1$ out of $D_j$; notice that $i\ddbar\rho_j$ extends to a \textit{smooth} real (1,1)-form on the whole $X$, the class of which is $2\pi c_1([D_j])$. 
We can also assume that $|\sigma_j|_j$ is constant near $D_k$ when $D_j\cap D_k =\varnothing$.
Let $\lambda$ be a nonnegative real parameter, and set $u_j:=\log(\lambda+\rho_j)=\log\big(\lambda-\log(|\sigma_j|^2_j)\big)$. 
Then for $\lambda$ big enough,
 \begin{equation}  \label{def_omega}
  \omega:=\omega_0-dd^c\mathfrak{u}
         =\omega_0-\sum_{j=1}^N dd^c u_j, \text{ where }\mathfrak{u}=\sum_{j=1}^N u_j=\sum_{j=1}^N\log\big(\lambda-\log(|\sigma_j|^2_j)\big) ,
 \end{equation}
defines a genuine Kähler form on $\XD$, that we will take as a reference metric in what follows ---we shall also assume that $\lambda=0$, since this is equivalent to replace the $|\cdot|_j$ by $e^{-\lambda}|\cdot|_j$. 
Indeed if $U$ is a polydisc of coordinates $(z_1,\dots, z_m)$ around some point of $D$ such that $U\cap D=\{z_1\cdots z_k=0\}$, then $\omega$ is mutually bounded near the divisor with $\sum_{j=1}^k\tfrac{idz_j\wedge d\overline{z_j}}{|z_j|^2\log^2(|z_j|^2)}+\sum_{j=k+1}^m idz_j\wedge d\overline{z_j}$, and moreover has bounded derivatives at any order with respect to this local model metric. 
Sharper asymptotics are described in Proposition 1.2 of \cite{auv}. 

Such a metric is complete and has finite volume on $\XD$, equal to $\tfrac{[\omega_0]^m}{m!}$. 
We are interested into generalizing the behaviour of our model $\omega$ and the reference potential $\mathfrak{u}$ to a whole class; for this we state:

 \begin{df}  \label{df_poinc_met}
  Let $\overline{\omega}$ be a locally smooth closed real (1,1) form on $\XD$. 
We say that $\overline{\omega}$ is a \textit{Kähler metric of Poincaré type in the class $\Omega=[\omega_0]_{dR}$}, denoted by $\overline{\omega}\in\mom$, if:
  \begin{itemize}  \setlength{\itemsep}{0pt}
   \item[(1)] $\overline{\omega}$ is quasi-isometric to $\omega$, i.e. $c\omega\leq \overline{\omega}\leq c^{-1}\omega$ on $\XD$ for some $c>0$; 
   \item[(2)] $\overline{\omega}=\omega+dd^c\psi$ for function $\psi\in \mathcal{E}$,
  \end{itemize}
where 
  \begin{equation*}  \label{eq_df_E}
    \mathcal{E}:=\big\{\varphi\in C^{\infty}_{loc}(\XD)|\,\varphi=O(\mathfrak{u})\,\, \text{and}\,\, |\nabla_{\omega}^j \varphi|\,\, \text{is bounded on $\XD$ for any } j\geq 1\big\}.
   \end{equation*}
Similarly, we denote by $\tmom$ the space of potentials of such metrics, computed with respect to $\omega$. 
 \end{df}

A fundamental representant of such a class is the Kähler-Einstein metric Tian and Yau \cite{tian-yau1} produce when $K[D]$ is ample and $\Omega=\mu c_1\big(K[D]\big)$, $\mu>0$.

It is obvious that provided $\psi\in\mathcal{E}$ such that $\|dd^c\psi\|_{\omega}$ is small enough and $\overline{\omega}\in\mom$, then $\overline{\omega}+dd^c\psi$ is again a metric lying in $\mom$. 
In this way $\tmom$ is nothing but an open neighbourhood of 0 in $\mathcal{E}$, and we will use in what follows the notation
 \begin{equation}
  \omega_{\psi}:=\omega+dd^c\psi \in\mom
 \end{equation}
for all $\psi\in\tmom$.
 
Besides their completeness and finite volume property, the metrics in $\mom$ have also in common that their Ricci forms lie in $-2\pi c_1\big(K[D]\big)$ (as an $L^2$ class). 
Consequently they all have the same mean scalar curvature, which we denote by:
 \begin{equation*}
  \overline{\scal}=-4\pi m\frac{c_1(K[D])[\omega_0]^{m-1}}{[\omega_0]^{m}}.
 \end{equation*}  

~

Finally, one observes when $D$ is smooth that for every $j\in\{1,\dots,N\}$, $\omega$ induces a Kähler form $\omega|_{D_j}$ on $D_j$. 
The class of this form is actually $[\omega_0|_{D_j}]$, which depends only on $\Omega$, and any Kähler form in this class has its Ricci form lying in $-2\pi c_1(K_{D_j})$, which can also be written as $-2\pi c_1\big(K[D]|_{D_j}\big)$, according to the adjunction formula and the triviality of $[D_k]$ over $D_j$ when $k\neq j$. 
Thus, any metric on $D_j$ in the class $[\omega_0|_{D_j}]$ gets as mean scalar curvature
 \begin{equation} \label{eq_mn_sc_crv}
  \overline{\scal}_{D_j}:=-4\pi(m-1)\frac{c_1(K[D]|_{D_j})[\omega_0|_{D_j}]^{m-2}}{[\omega_0|_{D_j}]^{m-1}}.
 \end{equation}

When $D$ is not smooth and $D_j$ intersects other components, the metric induced on $D_j\backslash \bigcup_{j'\neq j} D_{j'}$ is of Poincaré type, hence has its Ricci form in the class $-2\pi c_1\big(K_{D_j}[D_j\cap \sum_{j'\neq j}D_{j'}]\big)$, which writes again $-2\pi c_1\big(K[D]|_{D_j}\big)$, and thus the mean scalar curvature one gets is 
again given by \eqref{eq_mn_sc_crv}, which we denote again by $\overline{\scal}_{D_j}$. 
At last, we associate in both cases this number to $[\omega_0|_{D_j}]$ or $\mathcal{PM}_{[\omega_0|_{D_j}]}$ as \textit{the mean scalar curvature of this class}. 

We can now state the main result of this note:
 \begin{thm} \label{thm_topobstr}
  Let $\overline{\scal}_{D_j}$ the mean scalar curvature associated to the (Poincaré) class of $\omega|_{D_j}$ for $j=1,\dots, N$ as above. 
Suppose there exists a Kähler metric of Poincaré type in the class of $\omega_0$ (in the sense of Definition \ref{df_poinc_met}) with constant scalar curvature on $\XD$. 
Then for every $j$, $\overline{\scal}<\overline{\scal}_{D_j}$. 
In other words, we have for every $j$ the inequality
  \begin{equation*}
   m\frac{c_1\big(K[D]\big)[\omega_0]^{m-1}}{[\omega_0]^{m}} > (m-1)\frac{c_1\big([D_j]\big)c_1\big(K[D]\big)[\omega_0]^{m-2}}{c_1\big([D_j]\big)[\omega_0]^{m-1}}
  \end{equation*}
in terms of classes defined on $X$.
 \end{thm}

The present theorem genuinely deals with Poincaré type Kähler metrics specificity, since the constraint it states becomes empty in the compact case. 
Moreover, it represents a first step toward the notion of K-stability G. Székelyhidi \cite{sze} suggested for pairs $(X,D)$, formulated to take into account such a Poincaré behaviour. 
We shall explain these links in next part, first section. 
Nonetheless, the content of the theorem is rather intuitive, and can just be understood as a negative contribution in the scalar curvature from the component of the metric normal to the divisor, which looks like the Einstein metric of negative Ricci form on Poincaré's punctured disc. 

Let us give a few comments about the choice we made in Definition \ref{df_poinc_met} for the space of potentials. 
On the one hand, we asked for a priori $L^{\infty}$ control on the derivatives of order $\geq1$ because this comes into consideration in the analysis on potentials, especially in a description of $\XD$ near infinity we need to prove our theorem.
On the other hand, we allow non-sharp asymptotics near the divisor, to avoid the asymptotically product situation. 
Indeed, we shall see in part \ref{part_link}, second section, that the existence of an asymptotically product Kähler metric on $\XD$ implies that $D$ carries Kähler metrics with constant scalar curvature, and this implies the K-stability of $D$, which would not be clearly implied by the K-stability of $(X,D)$. 
Working with what Poincaré type metrics, on which sharp asymptotics are not required, thus allows us to avoid such an assumption on the divisor.  
Now in this framework Theorem \ref{thm_topobstr} is no more an immediate consequence of the hypothesis, as it would be if dealing with asymptotically product metrics; we also look at this point in part \ref{part_link}.

Finally, let us precise that Theorem \ref{thm_topobstr} is the consequence of growth constraints on the potential of 
a metric of Poincaré type described in Propositions \ref{prop_calcul1} and \ref{prop_calcul2} below, together with the use of the constant scalar curvature property. 
Since we need a suitable description of $X$ near $D$ to state those propositions, we shall deal with such a description in part \ref{part_fibr}, and state our propositions in part \ref{part_prop}, as well as we shall show how they imply Theorem \ref{thm_topobstr} and give their proofs. 
For simplicity, we assume that $D$ is smooth until the end of part \ref{part_prop}, and deal with the generalization to the simple normal crossings case in part \ref{part_gen}.

 \section{Links with Székelyhidi's suggestions}  \label{part_link}

   \subsection{K-stability for a triple $(X,D,L)$}

In \cite[\S 3.1.2]{sze}, G. Székelyhidi considers (a subclass of) Kähler metrics of Poincaré type resulting from a polarisation $L\to X$, the class of \textit{asymptotically hyperbolic metrics}. 
He then suggests the following conjecture, extending Yau-Tian-Donaldson of the compact case:
 \begin{conj}[Székelyhidi]
  Assume $D$ is smooth. 
  The triple $(X,D,L)$ is K-stable if, and only if, there exists an asymptotically hyperbolic Kähler metric in the polarization class.
 \end{conj}

Precisions on asymptotically hyperbolic metrics are given further.
Now, by contrast to the K-stability in the compact case, the K-stability for triples $(X,D,L)$ requires two types of conditions:
 \begin{enumerate}
  \item for any \textit{test-configuration} of $(X,D,L)$, the \textit{Futaki invariant} is non-negative, and vanishes only for product configurations;
  \item one can compute numbers $c_0\neq0$, $c_1$, $\alpha_1\neq0$ and $\alpha_2$ for any test-configuration of $(X,D,L)$; we ask for the inequality
   \begin{equation}  \label{ineq_cond2}
    \frac{c_1}{c_0}<\frac{\alpha_2}{\alpha_1},
   \end{equation}
 \end{enumerate}
Let us say that condition 1 is analogous to the computations involved in the compact case. 
On the other hand, condition 2 is original, and is meant to guarantee that the metrics in game do have a "Poincaré behaviour" near the divisor; in the following lines, we are looking deeper into this condition. 

Let us mention now we will not need a precise definition of test-configurations, and that we will not use Futaki invariants; the reader is referred to \cite{sze} for the details. 
We just need to know that $(X\times\C,D\times\C)$, above $(X,D)$, polarized by the pull-back of $L$ for the obvious projection and with the trivial $\C^*$-action, is a test-configuration for $(X,D,L)$ ; this is indeed the trivial one. 

In this way we give a few precisions on the numbers $c_0$, $c_1$, $\alpha_1$ and $\alpha_2$ computed for the trivial configuration. 
By definition, $c_0$ and $c_1$ are given by 
\begin{equation*}
  \begin{aligned}
   \frac{d_k+\tilde{d}_k}{2} =: c_0k^{m}+c_1 k^{m-1} + O(k^{m-2}), 
  \end{aligned}
 \end{equation*}
where $d_k$ (resp. $\tilde{d}_k$) is the dimension of $H^0(X, L^k)$ (resp. of $H^0\big(X, L^k_0\otimes\mathcal{O}(-D)\big)$). 

Similar computations are used to define $\alpha_1$ and $\alpha_2$, which verify the relation:
 \begin{equation*}
  \dim H^0(D,L^k_0|_{D_0})=:\alpha_1 k^{m-1}+\alpha_2 k^{m-2}+O(k^{m-3})
 \end{equation*}

There is no difficulty in computing those four terms for our trivial configuration:
 \begin{prop} \label{prop_cfg_tst_triviale}
  Consider the trivial test-configuration for the triple $(X,D,L)$. 
  Then:
   \begin{equation*}
    \begin{aligned}
     c_0      = &  \frac{c_1(L)^m}{m!}                ,&      &   c_1     = -\frac{c_1(L)^{m-1}\cdot c_1(K)+c_1(L|_D)^{m-1}}{2(m-1)!},   \\
     \alpha_1 = &  \frac{c_1(L|_D)^{m-1}}{(m-1)!}     ,&      &  \alpha_2 = -\frac{c_1(L|_D)^{m-2}\cdot c_1(K_D)}{2(m-2)!}   .
    \end{aligned}
   \end{equation*}
 \end{prop}

~

Now if one takes $[\omega_0]=2\pi c_1(L)$ in this case, one thus has $[\omega_0|_D]=2\pi c_1(L|_D)$, hence $[\omega_0|_D]^{m-1}=[\omega_0]^{m-1}\cdot c_1([D])$. 
A straightforward computation shows that the inequality $\tfrac{c_1}{c_0}<\tfrac{\alpha_2}{\alpha_1}$ corresponds to $m\tfrac{c_1(K[D])[\omega_0]^{m-1}}{[\omega_0]^{m}}>(m-1)\tfrac{c_1(K_{D})[\omega_0|_{D}]^{m-2}}{[\omega_0|_{D}]^{m-1}}$, which is the inequality stated in Theorem \ref{thm_topobstr} when $D$ is reduced to one component. 
Moreover, when $D$ is smooth and admits several components, the inequality of \eqref{ineq_cond2} can be obtained as the average of the inequalities of Theorem \ref{thm_topobstr}. 

To sum it up, Theorem \ref{thm_topobstr} is a first step toward the implication "existence of asymptotically hyperbolic Kähler metric $\Rightarrow$ K-stability of $(X,D,L)$".

~

\noindent\prf \textit{of Proposition \ref{prop_cfg_tst_triviale}.} 
According to Riemann-Roch theorem, for $k$ going to infinity, 
 \begin{equation*}
  d_k=\frac{c_1(L)^m}{m!}k^{m}-\frac{c_1(L)^{m-1}\cdot c_1(K)}{2(m-1)!} k^{m-1} + O(k^{m-2}),
 \end{equation*}
and 
 \begin{equation}  \label{eq_h0_Lk|D}
   h^0(D,L^k|_{D})=\frac{c_1(L|_D)^{m-1}}{(m-1)!}k^{m-1}-\frac{c_1(L|_D)^{m-2}\cdot c_1(K_D)}{2(m-2)!} k^{m-2} + O(k^{m-3}),
 \end{equation}
thus $\alpha_1=\tfrac{c_1(L|_D)^{m-1}}{(m-1)!}$ and $\alpha_2=-\tfrac{c_1(L|_D)^{m-2}\cdot c_1(K_D)}{2(m-2)!}$. 

There remains to compute $\tilde{d}_k$; now for $k$ big enough $H^1\big(X,L^k\otimes \mathcal{O}(-D)\big)=0$, and consequently the sequence 
 \begin{equation*} \label{eq_suite_exacte}
  0  \longrightarrow  H^0\big(X,L^k\otimes \mathcal{O}(-D)\big) 
     \longrightarrow  H^0(X,L^k)
     \longrightarrow  H^0(D,L^k|_D)
     \longrightarrow  0
 \end{equation*}
is exact, thus $\tilde{d}_k=d_k-h^0(D,L^k|_{D})$. 
Finally from \eqref{eq_h0_Lk|D} one has $c_0=\tfrac{c_1(L)^m}{m!}$, and $c_1=-\tfrac{1}{2}\big(\tfrac{c_1(L)^{m-1}\cdot c_1(K)}{(m-1)!}+\tfrac{c_1(L|_D)^{m-1}}{(m-1)!}\big)$. 
\cqfd

~

\subsection{Asymptotically hyperbolic Kähler metrics}

Our next comment on Székelyhidi's conjecture will concern more precisely the class of asymptotically hyperbolic near $D$ Kähler metrics. 
We start by the precise definition, which we take in \cite{sze}. 

Consider an open set $U$ of holomorphic coordinates $\{z_1,\dots,z_m\}$ around some point of $D$ (assumed smooth) such that $D=\{z_1=0\}$ in $U$. 
Then near $D$ an asymptotically hyperbolic metric $g$ looks like a product metric
 \begin{equation}  \label{eq_df_sz}
  \hat{g}_U=K\frac{|dz_1|^2}{|z_1|^2\log^2(|z_1|^2)}+h_U
 \end{equation}
at any order, i.e. $\big|\nabla^k_{\hat{g}_U}(g-\hat{g}_U)\big|_{\hat{g}_U}=o(1)$ near $D$ for all $k\geq0$. 
Here $K$ is a smooth positive function on $X$, and $h_U$ is a smooth extension of a metric on $D$. 

Clearly, asymptotically hyperbolic metrics are of Poincaré type, but they are actually rather specific, because their definition implies that:
 \begin{enumerate}
  \item $K$ is constant on $D$, and can hence be taken constant near $D$, the error being included in the $o(1)$ above;
  \item if $g$ has constant scalar curvature, so does the metric $g|_D$ it induces on $D$, and $\scal(g)<\scal(g|_D)$.
 \end{enumerate}
Fact 1 is proved by looking at the Kähler form $\omega_g$ of $g$; indeed, in the coordinates used above on $U$, for any $j,k\in\{2,\dots,m\},$
 \begin{align*}
  &\omega_{1\bar{1}} = \frac{K+o(1)}{|z_1|^2\log^2(|z_1|^2)},         &                &\omega_{1\bar{k}}= O\big(\frac{1}{|z_1||\log|z_1||}\big),  \\  
  &\omega_{j1}      = O\big(\frac{1}{|z_1||\log|z_1||}\big),          &\text{and}\quad &\omega_{j\bar{k}}= (\omega_{h_U})_{j\bar{k}}                
 \end{align*}
at any order, in asymptotically hyperbolic sense (or Poincaré, because of mutual bounds). 
Now the Kähler hypothesis tells us that for  any $j,k\in\{2,\dots,m\}$, $\partial_{z_j}\omega_{1\bar{1}}=\partial_{z_1}\omega_{j\bar{1}}=O\big(\tfrac{1}{|z_1||\log|z_1||}\big)$. 
Now $\partial_{z_j}\omega_{1\bar{1}}=\frac{\partial_{z_j}K+o(1)}{|z_1|^2\log^2(|z_1|^2)}$, and hence ($K$ is smooth on $X$), $\partial_{z_j}K\equiv0$ on $D\cap U$. 
Thus $K$ is constant on $D\cap U$; since such functions $K$ may depend on the open set of coordinates but have to patch, we deduce that they are constant on any component of $D$ with common value, at least in formula \eqref{eq_df_sz}. 

All this can be reformulated saying that on any $U$ as above, $\omega_g=\tfrac{Aidz_1\wedge d\overline{z_1}}{|z_1|^2\log^2(|z_1|^2)}+\omega_{g|_D}+o(1)$ for some positive constant $A$ independent of $U$, the perturbation being understood at any order in asymptotically hyperbolic metric.

Fact 2 readily follows from this latter observation. 
Indeed, one has $\omega_g^m=m\tfrac{Aidz_1\wedge d\overline{z_1}}{|z_1|^2\log^2(|z_1|^2)}\wedge(\omega_g|_D)^{m-1}$.
The Ricci form $\varrho_g$ of $g$ is thus given by 
 \begin{equation*}
  \varrho_{g|_D}-i\ddbar\log\Big(\frac{A}{|z_1|^2\log^2(|z_1|^2)}\Big)+o(1)=\varrho_{g|_D}-\frac{2idz_1\wedge d\overline{z_1}}{|z_1|^2\log^2(|z_1|^2)}+o(1).
 \end{equation*} 
Now $\scal(g)\omega_g^m=2m\varrho_g\wedge\omega_g^{m-1}$, which develops into
 \begin{align*}
  \scal(g)\omega_g^m& = m\Big[2(m-1)\varrho_{g|D}\wedge\omega_{g|D}^{m-2}-4A^{-1}\omega_{g|D}^{m-1}\Big]\wedge\frac{Aidz_1\wedge d\overline{z_1}}{|z_1|^2\log^2(|z_1|^2)}+o(1)\\
                    & = m\big(\scal(g|_D)-4A^{-1}\big)\omega_{g|_D}^{m-1}\wedge\frac{Aidz_1\wedge d\overline{z_1}}{|z_1|^2\log^2(|z_1|^2)}+o(1)
 \end{align*}
i.e. $\scal(g)=\scal(g|_D)-4A^{-1}+o(1)$. 
Thus if $\scal(g)$ is constant, $\scal(g|_D)$ is too, the $o(1)$ drops, $\scal(g|_D)>\scal(g)$ and $A=\tfrac{4}{\scal(g|_D)-\scal(g)}$.

~

These computations illustrate the intuitive interpretation of Theorem \ref{thm_topobstr}. 
Of course, such computations are no longer possible in the absence of asymptotics like those coming from the definition of asymptotically hyperbolic metrics, and this is why we develop hereafter some techniques to get our result in the more general class of Poincaré type Kähler metrics.

\section{A fibration near the divisor}  \label{part_fibr}

To begin with, let us suppose that $D$ is reduced to one (smooth) component, and denote by $\sigma$ a defining section for $D$. 

We consider a tubular neighbourhood $\mathcal{N}_A$ of $D$ ($A$ is a real parameter to be fixed), with projection $p$, obtained from the exponential map associated to a smooth metric on $X$, $\omega_0$ say. 
On $\mathcal{N}_A$, an $S^1$ action comes form the identification of $\mathcal{N}_A$ with some neighbourhood $\mathcal{V}$ of the null section of the normal holomorphic bundle $N_D=\tfrac{T^{1,0}X|_D}{T^{1,0}D}$ and this action leaves invariant the projection $p:\mathcal{N}_A\simeq\mathcal{V}\subset N_D\to D$. 
Now we complete $p$ by making the function $\mathfrak{u}$, which is reduced to $\log\big(-\log(|\sigma|^2)\big)$, invariant under the circle action. 
To be more precise, let $T$ the infinitesimal generator of the action, with flow $\Phi_{\vartheta}$. 
We set: 
 \begin{equation*}
  t:=\log\Big[\log\Big(-\frac{1}{2\pi}\int_{S^1} \Phi_{\vartheta}^*(|\sigma|^2)\,d\vartheta\Big)\Big]
 \end{equation*}
near $D$, and extend it smoothly away from the divisor. 
If we denote the couple $(p,t)$ by $q$, we have the following diagram:
 \begin{equation} \label{eqfibr}
  \xymatrix{
    S^1 \ar[r] & \mathcal{N}_A\backslash D \ar[d]^{q=(t,p)} \\
               & [A,+\infty)\times D
  }
 \end{equation}
It is easy to see that $t=\mathfrak{u}$ up to a perturbation which is a $O(e^{-t})$, that is, a $O\big(\tfrac{1}{|\log|\sigma||}\big)$, as well as its derivatives of any order (in Poincaré type metric).
Finally, $A$ and $\mathcal{N}_A$ are adjusted so that $\mathcal{N}_A\backslash D=\{t\geq A\}\subset\XD$.

We associate to the circle action on $\mathcal{N}_A$ a connection 1-form $\eta$, as follows: if $g$ is of Poincaré type (e.g. $g$ is the Riemannian metric associated to $\omega$) and keeping $T$ as the infinitesimal generator of the action with flow $\Phi_{\vartheta}$, we set at any point $x$ of $\mathcal{N}_A$
 \begin{equation*}
  \hat{\eta}_x= \int_0^{2\pi}\Phi_{\vartheta}^*\Big(\frac{g_x(\cdot,T)}{g_x(T,T)}\Big) d\vartheta \quad \text{ and } \quad \eta_x=2\pi\Big(\int_{S^1}\hat{\eta}\Big)^{-1}\hat{\eta}_x,
 \end{equation*}
where the $S^1$ of the last integral is the fiber $q^{-1}(x)$. 
In this way, for all $x\in\mathcal{N}_A$, $\int_{q^{-1}(x)}\eta=2\pi$. 

Moreover, if we consider around some point of $D$ a neighbourhood of holomorphic coordinates $(z_1=re^{i\theta},\dots,z_m)$ such that $D$ is given by $z_1=0$, one has $\eta=d\theta$ up to a term which is a $O(1)$ as well as its derivatives of any order with respect to $\omega$. 
Hence if we assume that $g$ denotes the Riemannian metric associated to $\omega$, one has
 \begin{equation}  \label{eq_asmpt_g}
  g=dt^2+4e^{-2t}\eta^2+p^{*}g_D+O\big(e^{-t}\big)
 \end{equation}
with $g_D$ the metric associated to $\omega_0|_D$, where $O\big(e^{-t}\big)$ is understood at any order in Poincaré metric; this follows from \cite[Proposition 1.2]{auv}. 
This means for example that $Jdt=2e^{-t}\eta+O\big(e^{-t}\big)$, this $O\big(e^{-t}\big)$ being understood in the same way. 

One can furthermore use the fibration \eqref{eqfibr} as follows. 
Let $f\in C^{k,\alpha}\big(\XD\big)$ (in the $\alpha=0$ case, which is relevant here, $C^{k,\alpha}\big(\XD\big)$ is defined with help of $\nabla_{\omega}$; see \cite[section 1.2]{auv} for the definition of such Hölder spaces when $\alpha\in(0,1)$ ); one writes the decompositions
 \begin{equation}  \label{eqdecomp1}
  f=(\Pi_0 f)(t,z)+\Pi_{\perp}f = f_0(t)+f_1(t,z)+\Pi_{\perp}f
 \end{equation}
where $z=p(x)$, with:
 \begin{equation*}
  (\Pi_0 f)(t,z)=\frac{1}{2\pi}\int_{q^{-1}(x)}f\,\eta, \quad \text{and} \quad f_0(t)=\frac{1}{\vl(D)}\int_D (\Pi_0 f)(t,z)\,\vol^{g_D},
 \end{equation*}
and $\vl(D)$ computed with respect to $g_D$, hence equal to $\tfrac{[\omega_0|_D]^{m-1}}{(m-1)!}$, or $\tfrac{c_1([D])\cdot[\omega_0|_D]^{m-1}}{(m-1)!}$ by Lelong's formula. 

Using equation \eqref{eqdecomp1} and the definition of $C^{k,\alpha}\big(\XD\big)$, since the $S^1$ fibers are of length equivalent to $e^{-t}$ for $g$, it is not difficult to see that on an open set of holomorphic coordinates as above, if $j\leq k$,
 \begin{equation*}
  \mathcal{D}_{\ell,j-\ell}\big(\Pi_{\perp}f\big) =O(e^{-(k-\ell+\alpha)t})
 \end{equation*}
as soon $\mathcal{D}_{\ell,j-\ell}$ denotes a product $(j-\ell)$ factors of which are equal to $e^t\partial_{\theta}$, and the remaining $\ell$ factors are in $\{r|\log r|\partial_r,\partial_{z_\beta},\dbar_{z_\beta}, \,\beta\geq2\}$, were $r=|z_1|$.

If moreover $J_D$ is the complex structure on $D$, we have that $p^*J_D$ differs of $J$ restricted to $\bigoplus_{\beta\geq2} \big(\C\tfrac{\partial}{\partial z_{\beta}}\oplus\C\tfrac{\partial}{\partial \overline{z_{\beta}}}\big)$ by a perturbation which is a $O(e^{-t})$ as well as its derivatives at any order. 
An application of these estimations is that one has, if $f\in C^{k,\alpha}\big(\XD\big)$, $k\geq 2$, $\alpha\in(0,1)$,
 \begin{equation*} \label{eq_deriv1}
  df   = \big(\partial_t f_0(t)+\partial_t f_1(t,z)\big) dt+d_Df_1(t,z)+O(e^{-(k-1+\alpha)t})
 \end{equation*}
 \begin{equation} \label{eq_deriv2}
  d^cf = Jdf = 2\big(\partial_t f_0(t)+\partial_t f_1(t,z)\big)e^{-t}\eta+d_D^cf_1(t,z)+O(e^{-t})
 \end{equation}
with $d_D^c=p^*(J_Dd)$, and
 \begin{equation}  \label{eq_deriv3}
  \begin{aligned}
   dd^cf  = & 2\big(\partial_t^2 f_0(t)+\partial_t^2 f_1(t,z)-\partial_t f_0(t)-\partial_t f_1(t,z)\big)e^{-t}dt\wedge\eta \\
            & + 2e^{-t}d_D\big(\partial_t f_0(t)+\partial_t f_1(t,z)\big)\wedge\eta +dt\wedge d_D^c\partial_tf_1(t,z) \\
            & + dd^c_Df_1(t,z)+O\big(e^{-t}\big)
  \end{aligned}
 \end{equation}
with $dd_D^c=p^*(dJ_Dd)$

Eventually, if we replace $f\in C^{k,\alpha}\big(\XD\big)$ by a potential $\varphi$ of a metric in $\mom$, computed with respect to $\omega$, decomposition \eqref{eqdecomp1} and estimates \eqref{eq_deriv2} and \eqref{eq_deriv3} apply again, except that $\varphi_0(t)=O(t)$, and $\partial_t^j\varphi_0=O(1)$ for all $j\geq 1$. 

~

All this description is easily transposable in the case the number $N$ of disjoint components of $D$ is strictly bigger than $1$, by working near one component and away from the others; 
in this case we shall add to the objects above an index $j$ (e.g. $t_j$, $\eta_j$, $\varphi_{0,j}$, etc.) to specify which component of the divisor they refer to. 

~

 \section{The key propositions} \label{part_prop}

 \subsection{Statements, and proof of Theorem \ref{thm_topobstr} (smooth divisor case)}

We assume again in this part that $D$ is smooth.
We come now to the statement of the two key propositions of which Theorem \ref{thm_topobstr} is a corollary: assuming $\omega_{\varphi}:=\omega+dd^c\varphi$ has constant scalar curvature (hence equal to $\overline{\scal}=-4\pi m\tfrac{c_1(K[D])[\omega_0]^{m-1}}{[\omega_0]^{m}}$ at any point), one find constraints on $\varphi$ which will translate into constraints on $\overline{\scal}$. 
For this, we state first:
 \begin{prop}[$D$ smooth]  \label{prop_calcul1}
  Assume $\varphi\in\tmom$. 
Fix $j\in\{1,\dots,N\}$, and consider the compact subdomains $\{t_j\leq s\}\subset X\backslash D_j$. 
Then $\bigcup_{s\geq0}\{t_j\leq s\} = X\backslash D_j$, the integral $\int_{\XD}e^{t_j}\omega_{\varphi}^m$ diverges to $+\infty$, and: 
   \begin{equation}  \label{eq_calcul1}
    \int_{\{t_j\leq s\}}e^{t_j}\omega_{\varphi}^m = 4\pi m!\vl(D_j)\big(s-\varphi_{0,j}(s)\big)+O(1)
   \end{equation}
 when $s$ goes to $\infty$, for $j=1,\dots,N$. 
 \end{prop}

~

We also state:
  \begin{prop}[$D$ smooth]  \label{prop_calcul2}
  Assume $\varphi\in\tmom$. 
Then: 
   \begin{equation}  \label{eq_calcul2}
 \int_{\{t_j\leq s\}}\scal_{\varphi}e^{t_j} \omega_{\varphi}^m  =4\pi m!\vl(D_j)\big((\overline{\scal}_{D_j}-2)s-\overline{\scal}_{D_j}\varphi_{0,j}(s)\big)+O(1)
   \end{equation}
 when $s$ goes to $\infty$, for $j=1,\dots,N$. 
 \end{prop}

The proofs are postponed to section \ref{section_prfs}, since we shall see for now how Propositions \ref{prop_calcul1} and \ref{prop_calcul2} readily imply Theorem \ref{thm_topobstr} when $D$ is smooth.

~

\noindent\prf~\textit{of Theorem \ref{thm_topobstr} from Propositions \ref{prop_calcul1} and \ref{prop_calcul2}}. 
Assume $\omega_{\varphi}=\omega+dd^c\varphi$ has constant scalar curvature. 
Fix $j\in\{1,\dots,N\}$.
Compare equations \eqref{eq_calcul1} and \eqref{eq_calcul2}, after having multiplied the first one by $\overline{\scal}$ (which equals $\scal_{\varphi}$ at every point of $\XD$). 
It follows, for $s$ going to $\infty$:
 \begin{equation*}
  (\overline{\scal}_j-\overline{\scal}-2)s = (\overline{\scal}_{j}-\overline{\scal})\varphi_{0,j}(s)+O(1),
 \end{equation*}
or $(\overline{\scal}_j-\overline{\scal})\big(s-\varphi_{0,j}(s)\big)=2s+O(1)$ to be more explicit (here and from now on, $\overline{\scal}_j$ stands for $\overline{\scal}_{D_j}$). 
Having $\overline{\scal}_j=\overline{\scal}$, \textit{i.e.} $2s=O(1)$, would thus be absurd. 
It follows that $\overline{\scal}\neq\overline{\scal}_j$, and $s-\varphi_{0,j}(s)= \tfrac{2}{\overline{\scal}_j-\overline{\scal}}s+O(1)$. 
Now, in order to determine the sign of $\overline{\scal}_{j}-\overline{\scal}$, we use the divergence of $\int_{\XD}e^{t_j}\omega_{\varphi}^m$.
Since $\int_{\XD}e^{t_j}\omega_{\varphi}^m$ is the increasing limit of $\int_{\{t_j\leq s\}}e^{t_j}\omega_{\varphi}^m$, from the asymptotics \eqref{eq_calcul1}, we know that $\big(s-\varphi_{0,j}(s)\big)$ tends to $+\infty$ when $s$ goes to $\infty$. 
This would not be compatible with an inequality $\overline{\scal}>\overline{\scal_j}$, hence $\overline{\scal}<\overline{\scal}_j$. 

We thus have proved the inequality $m\tfrac{c_1(K[D])[\omega_0]^{m-1}}{[\omega_0]^{m}}>(m-1)\tfrac{c_1(K_{D_j})[\omega_0|_{D_j}]^{m-2}}{[\omega_0|_{D_j}]^{m-1}}$. 
To recover the inequality of the theorem, one uses that given $(m-1)$ classes $[\alpha_1],\dots,[\alpha_{m-1}]$ with (1,1) representatives $\alpha_1,\dots,\alpha_{m-1}$ on $X$, then the cup product $[\alpha_1]\cdots[\alpha_{m-1}]c_1([D_j])$ is equal to $[\alpha_1|_{D_j}]\cdots[\alpha_{m-1}|_{D_j}]$.
One also uses the adjunction formula $K[D]|_{D_j}\approx K[D_j]|_{D_j} \approx K_{D_j}$ (the $[D_k]$, $k\neq j$, are trivial over $D_j$).
\hfill $\square$

 \begin{rmk}
  This proof also gives us a sharper description of the potential $\varphi$ of a Poincaré metric with constant scalar curvature; indeed, in a nutshell, one has $\varphi-\sum_{j=1}^N a_jt_j\in C^{\infty}(\XD)$, where $a_j=\tfrac{\overline{\scal_j}-\overline{\scal}-2}{\overline{\scal_j}-\overline{\scal}}$ for each $j$. 
In other words, we directly pass from "$\log\log$" terms (the $t_j$) to bounded terms in the development of such a potential near the divisor, whereas terms like $t_j^{1/2}$ are a priori authorized to contribute to potentials of Poincaré type. 
 \end{rmk}

~

\subsection{Proofs of Propositions \ref{prop_calcul1} and \ref{prop_calcul2}}  \label{section_prfs}

\subsubsection{Proof of Proposition \ref{prop_calcul1}} \label{pgrph_prf_prop1}
Fix $j\in\{1,\dots,N\}$.
The exhaustion $\bigcup_{s\geq0}\{t_j\leq s\} = X\backslash D_j$ is clear, arising from the construction of $t_j$, and more precisely from the estimate $t_j=u_j+o(1)$ where we recall that $u_j=\log\big(-\log(|\sigma_j|^2)\big)$ and $D_j=\{\sigma_j=0\}$ in $X$. 

Let us look at the divergence of $\int_{\XD}e^{t_j}\omega_{\varphi}^m$. 
As $e^{t_j}\omega_{\varphi}^m$ is mutually bounded with the volume form of the model $\omega$, which is $\tfrac{1}{m!}\omega^m$, near each $D_k$ different from $D_j$, and since $\omega$ has a finite volume on $\XD$, the integration near those $D_k$ does not contribute to the divergence. 
On the other hand, near $D_j$, $e^{t_j}\omega_{\varphi}^m$ is mutually bounded with $e^{t_j}\omega^m$, itself behaving as the cylindrical volume form $dt_j\wedge\eta_j\wedge p_j^*(\omega|_{D_j})^{m-1}$, of infinite volume. 

We now prove formula \eqref{eq_calcul1}. 
Set $\Theta=\omega^{m-1}+\omega^{m-2}\wedge\omega_{\varphi}+\cdots+\omega_{\varphi}^{m-1}$, so that $\omega_{\varphi}^m=\omega^{m-1}+dd^c\varphi\wedge\Theta$, and that $\Theta$ is a closed $(m-1,m-1)$-form.  
Our $j\in\{1,\dots,N\}$ is fixed again. 
By Stokes' theorem, we have for $s$ going to $\infty$ (assume $s$ is big enough so that $D_k\subset\{t_j\leq s\}$ when $k\neq j$):
 \begin{equation*}
  \int_{\{t_j\leq s\}} e^{t_j} \omega_{\varphi}^m = \int_{\{t_j\leq s\}} e^{t_j} \omega^m + \int_{\{t_j=s\}}e^{t_j}d^c\varphi\wedge\Theta- \int_{\{t_j\leq s\}}d(e^{t_j})\wedge d^c\varphi\wedge\Theta.
 \end{equation*}
Let us simplify the notation and drop the $j$ indexes; this will not lead to some confusion, since we are precisely dealing with integrals near $D_j$. 
Now as $\Theta$ has type $(m-1,m-1)$, $d(e^t)\wedge d^c\varphi\wedge\Theta= d\varphi\wedge d^c(e^t)\wedge\Theta$, and again by Stokes,
 \begin{equation*}
  \int_{\{t\leq s\}}d\varphi\wedge d^c(e^t)\wedge\Theta = \int_{\{t=s\}}\varphi d^c(e^t)\wedge\Theta- \int_{\{t\leq s\}} \varphi dd^c(e^t)\wedge\Theta,
 \end{equation*}
hence:
 \begin{equation} \label{eq_som_intgrl1}
  \int_{\{t\leq s\}} e^t \omega_{\varphi}^m = \int_{\{t\leq s\}} e^t \omega^m + \int_{\{t\leq s\}} \varphi dd^c(e^t)\wedge\Theta+ e^s\int_{\{t=s\}}\big(d^c\varphi-\varphi d^ct\big)\wedge\Theta
 \end{equation}

To get to \eqref{eq_calcul1}, we shall analyze the different summands in game. 
First, $e^t=\rho_j+O(1)$ at any order, hence $dd^c(e^t)= dd^c\rho_j+O(1)$; now, $dd^c\rho_j$ is bounded (for $\omega_0$, and thus for $\omega$, of Poincaré type) since it extends smoothly through $D$. 
In this way, as $\Theta$ is itself dominated by $\omega^{m-1}$ and $\varphi$ is $L^1$ for $\omega$, we deduce that $\int_{\{t\leq s\}} \varphi dd^c(e^t)\wedge\Theta$ converges when $s$ goes to infinity, and in particular is bounded. 

Let us pass to $\int_{\{t=s\}} d^c\varphi\wedge\Theta$. 
According to derivation formulas \eqref{eq_deriv2} and \eqref{eq_deriv3}, since $dt=0$ on $\{t=s\}$, one has on this slice
 \begin{equation*}
  \left\{ 
   \begin{aligned}
    \omega           & = p^*\omega_D+O(e^{-s}) \\
    \omega_{\varphi} & = p^*(\omega_{D_j}+ dd^c_{D_j}\varphi_1) + 2d_{D_j}\dot{\varphi_1}\wedge e^{-s}\eta+ O(e^{-s}) \\
    d^c\varphi       & = 2(\dot{\varphi_0}+\dot{\varphi_1})e^{-s}\eta+d^c_{D_j}\varphi_1 + O(e^{-s})
   \end{aligned}
  \right.
 \end{equation*}
where $\dot{}$ stands for $\partial_t$, the $O$ being still understood for Poincaré metrics. 
The notation $d_{D_j}\dot{\varphi_1}$ means $d(\dot{\varphi_1}|_{t=s})$, or the pull-back by $p$ of the differential of the function induced on $D_j$ by $\dot{\varphi_1}$ when fixing $t=s$; similarly $d^c_{D_j}\varphi_1$ stands for the pull-back of $J_{D_j}d (\varphi_1|_{t=s})$. 

In this way, we obtain on  $\{t=s\}$
 \begin{equation*} 
  \begin{aligned}
   d^c\varphi\wedge \Theta = & 2(\dot{\varphi_0}+\dot{\varphi_1})e^{-s}\eta\wedge p^*\big(\omega_{D_j}^{m-1}+\cdots+(\omega_{D_j}+idd^c\varphi_1)^{m-1}\big)\\
                             & -2\sum_{k=0}^{m-1} k e^{-s}\eta\wedge d_{D_j}\varphi_1\wedge d_{D_j}^c\varphi_1\wedge p^*\big((\omega_{D_j}+dd^c_{D_j}\varphi_1)^{m-1-k}\wedge\omega_{D_j}^{k-1}\big)\\
                             & +O(e^{-s}).
  \end{aligned}
 \end{equation*}
We shall notice that the latter $O(e^{-s})$ is a $(2m-1)$-form on $\{t=s\}$; we could for example write it as $\vareps(x)e^{-s}\eta\wedge (p^*\omega_{D_j})^{m-1}$, with $\vareps(x)=O\big(e^{-t(x)}\big)$. 

Up to some $O(e^{-s})$, $d^c\varphi\wedge \Theta$ is $S^1$-invariant, thus as the fibers are of length $2\pi$ for $\eta$, one has
 \begin{equation*}
  \begin{aligned}
   \int_{\{t=s\}}d^c\varphi\wedge \Theta =  4\pi e^{-s}\Big(& \int_{D_j} (\dot{\varphi_0}+\dot{\varphi_1})(s,\cdot)\big(\omega_{D_j}^{m-1}+\cdots+(\omega_{D_j}+idd^c\varphi_1)^{m-1}\big) \\
                                                            & - \int_{D_j} d\varphi_1(s,\cdot)\wedge d^c\varphi_1(s,\cdot)\wedge \Xi\big(\varphi_1(s,\cdot)\big) +O(e^{-s})\Big) 
  \end{aligned}
 \end{equation*}
where $\Xi(u)= \sum_{k=1}^{m-1} k \big((\omega_{D_j}+dd^c u)^{m-1-k}\wedge\omega_{D_j}^{k-1}\big)$ for any (smooth) function $u$ on $D_j$. 

In view of:
 \begin{itemize}
  \item the dependence on $s$ only of $\varphi_0$ ;
  \item the equality
   \begin{equation*}
    \int_{D_j} (\omega_{D_j}+dd^c_{D_j}\varphi_1)^{m-1-k}\wedge\omega_{D_j}^{k}=[\omega_{D_j}]^{m-1}=:(m-1)!\vl(D_j)
   \end{equation*}
   for all $k\in\{0,\dots,m-1\}$ ;
  \item the boundedness of the component $\varphi_1$ and its differential;
 \end{itemize}
it follows that: 
 \begin{equation}  \label{eq_int11}
  \int_{\{t=s\}}d^c\varphi\wedge \Theta = 4\pi e^{-s}m!\vl(D_j)\dot{\varphi_0}(s)+O(e^{-s}), 
 \end{equation}
that is: $e^s\int_{\{t=s\}}d^c\varphi\wedge \Theta =  4\pi\vl(D_j)m!\dot{\varphi_0}(s)+O(1)=O(1)$. 

Using a similar process, and setting $\Upsilon(u)=\omega_{D_j}^{m-1}+\cdots+(\omega_{D_j}+dd^c u)^{m-1}$ for $u$ any function on $D_j$, one gets : 
 \begin{equation} \label{eq_int12}
  \begin{aligned}
   e^s\int_{\{t=s\}} \varphi  d^ct\wedge\Theta 
            = & 4\pi m!\vl(D_j)\varphi_0(s) 
              + 4\pi \int_{D_j} \varphi_1(s,\cdot)\Upsilon\big(\varphi_1(s,\cdot)\big)+ O(e^{-s}) \\
            = & 4\pi m!\vl(D_j)\varphi_0(s) +O(1).
  \end{aligned}
 \end{equation}

There remains to analyze $\int_{\{t\leq s\}} e^t \omega^m $. 
We can write, in a neighbourhood of $D_j$, $\omega^m = 2me^{-t}dt\wedge \eta\wedge (p^*\omega_{D_j})^{m-1}\big(1+O(e^{-t})\big)$, according to \eqref{eq_asmpt_g} (or rather its analogue near $D_j$ in the case when $D$ has several components). 
This we write again $e^{t}\omega^m = 2mdt\wedge \eta\wedge (p^*\omega_{D_j})^{m-1}+O(e^{-t})$, $O(e^{-t})$ understood in the sense of volume forms of metrics of Poincaré type. 
One more use of the $S^1$-invariance of $dt\wedge (p^*\omega_{D_j})^{m-1}+O(e^{-t})$ on $\{t \geq A\}$ and the equality $\int_{D_j} \omega_{D_j}^{m-1}=(m-1)!\vl(D_j)$, we thus have for any $s$ big enough:
 \begin{equation} \label{eq_int13}
  \int_{\{t\leq s\}} e^t \omega^m = \int_{\{t\leq A\}} e^t \omega^m + \int_A^s\big(4\pi m!\vl(D_j)+O(e^{-t})\big)dt = 4\pi m!\vl(D_j)s+O(1).
 \end{equation}

The proposition is proved by collecting \eqref{eq_som_intgrl1}, \eqref{eq_int11}, \eqref{eq_int12} et \eqref{eq_int13}. \cqfd

\subsubsection{Proof of Proposition \ref{prop_calcul2}}  \label{pgrph_prf_prop2}

The techniques we use to get to estimate \eqref{eq_calcul2} are exactly similar to those we just used to prove Proposition \ref{prop_calcul2}.
However the starting point is the following formulas for the computation of the scalar curvature: $\scal_{\varphi}=2(\Lambda_{\varphi}\varrho_{\varphi})$, or equivalently $\scal_{\varphi}\omega_{\varphi}^m=2m\varrho_{\varphi}\wedge\omega_{\varphi}^{m-1}$ ;
here we denote by $\varrho_{\varphi}$ the Ricci from of $\omega_{\varphi}$, and we denote by $\varrho$ that of $\omega$. 
One has between those two forms the relation: $\varrho_{\varphi}= \varrho-\tfrac{1}{2}dd^c f$, where we set $f=\log\big(\tfrac{\omega_{\varphi}^m}{\omega^m}\big)\in C^{\infty}(\XD)$. 
Multiplying by $2m\omega_{\varphi}^{m-1}$ yields:
 \begin{equation*}
  \scal_{\varphi}\omega_{\varphi}^m = 2m \varrho\wedge\omega_{\varphi}^{m-1}-mdd^cf\wedge\omega_{\varphi}^{m-1}.
 \end{equation*}
that is: $\int_{\{t\leq s\}}e^t\scal_{\varphi}\omega_{\varphi}^m=2m\int_{\{t\leq s\}}e^t\varrho\wedge\omega_{\varphi}^{m-1}-m\int_{\{t\leq s\}}e^tdd^cf\wedge\omega_{\varphi}^{m-1}$ after integration (here again we fix $j\in\{1,\dots,N\}$, and drop it as an index when there is no risk of confusion).
Here one can write $\omega_{\varphi}^{m-1}=\omega^{m-1}+dd^c\varphi\wedge\Psi$, with $\Psi$ the $(m-2,m-2)$-form $\omega^{m-2}+\cdots+\omega_{\varphi}^{m-2}$. 
This gives the sum
 \begin{equation} \label{eq_som_intgrl2}
  \int_{\{t\leq s\}}e^t\varrho\wedge\omega_{\varphi}^{m-1} = \int_{\{t\leq s\}}e^t\varrho\wedge\omega^{m-1} + \int_{\{t\leq s\}}e^t\varrho\wedge dd^c\varphi\wedge\Psi.
 \end{equation}

In order to estimate $\int_{\{t\leq s\}}e^t\varrho\wedge\omega^{m-1}$, we notice that, $\varrho$ being asymptotically a product, its Ricci form will be too. 
More precisely, denoting by $\varrho_j$ the Ricci form of $\omega_{D_j}$, since from \eqref{eq_asmpt_g} one has $\omega=2dt\wedge e^{-t}\eta+p^*(\omega|_{D_j})+O(e^{-t})$, we get the asymptotics
 \begin{equation}  \label{eq_asmpt_ric}
  \varrho = p^*\varrho_j - dt\wedge e^{-t}\eta + O(e^{-t}), 
 \end{equation}
as sketched in our discussion in section \ref{part_link}, with less precise asymptotics. 
Hence near $D_j$:
 \begin{equation*}
  \varrho\wedge\omega^{m-1} = 2(m-1)dt\wedge e^{-t}\eta\wedge p^*\big(\varrho_j\wedge\omega_{D_j}^{m-2}\big)- 2dt\wedge e^{-t}\eta\wedge p^*\big(\omega_{D_j}^{m-1}\big) + O(e^{-t}).
 \end{equation*}
Proceeding as for $\int_{\{t\leq s\}} e^t \omega^m $, one gets this way:
 \begin{equation*}
  \int_{\{t\leq s\}}e^t\varrho\wedge\omega^{m-1} = 4\pi s\big((m-1)[\varrho_j][\omega_{D_j}]^{m-2}-[\omega_{D_j}]^{m-1}\big)+O(1),
 \end{equation*}
that is to say, since $[\varrho_{j}][\omega_{D_j}]^{m-2}=\tfrac{\overline{\scal}_{j}[\omega_{D_j}]^{m-1}}{2(m-1)}$ and $[\omega_{D_j}]^{m-1}=(m-1)!\vl(D_j)$,
 \begin{equation*} \label{eq_int21}
  2m\int_{\{t\leq s\}}e^t\varrho\wedge\omega^{m-1} = 4\pi m!s (\overline{\scal}_{j}-2)\vl(D_j)+O(1).
 \end{equation*}

We proceed by successive integrations by parts for the summand $\int_{\{t\leq s\}}e^t\varrho\wedge dd^c\varphi\wedge\Psi$ in equation \eqref{eq_som_intgrl2}:
 \begin{equation*}
  \int_{\{t\leq s\}}e^t\varrho\wedge dd^c\varphi\wedge\Psi = \int_{\{t\leq s\}}\varphi dd^c(e^t)\wedge\varrho \wedge\Psi+e^s\int_{\{t= s\}}(d^c\varphi-\varphi d^ct)\wedge\varrho \wedge\Psi.
 \end{equation*}

The first integral of the right-hand-side member is bounded for the same reasons than $\int_{\{t\leq s\}}\varphi dd^c(e^t) \wedge\Theta$ (see the proof of Proposition \ref{prop_calcul1}, \S \ref{pgrph_prf_prop1}). 
Using the asymptotics of $\varrho$ (formula \eqref{eq_asmpt_ric}) and the same arguments as for $e^s\int_{\{t= s\}}(\varphi d^ct-d^c\varphi)\wedge\Theta$ leads us to: 
 \begin{align*}
  e^s\int_{\{t= s\}}(\varphi d^ct-d^c\varphi)\wedge\varrho \wedge\Psi = & 2\pi (m-1)!\overline{\scal}_{j}\big(\dot{\varphi_0}(s)-\varphi_0(s)\big)+O(1) \\
                                                                      = & -2\pi (m-1)!\overline{\scal}_{j}\varphi_0(s)+O(1),
 \end{align*}
i.e. $2m\int_{\{t\leq s\}}\varphi dd^c(e^t)\wedge\varrho \wedge\Psi=-4\pi m!\overline{\scal}_{j}\varphi_0(s)+O(1)$

To conclude one has to see that $\int_{\{t\leq s\}}e^tdd^c f\wedge \omega^{m-1}$ remains bounded; here again two successive integrations by parts give:
 \begin{equation*}
  \int_{\{t\leq s\}}e^tdd^c f\wedge \omega^{m-1} = \int_{\{t\leq s\}}f dd^c(e^t) \wedge \omega^{m-1}+ e^s\int_{\{t= s\}}(d^cf-fd^ct)\wedge \omega^{m-1} ;
 \end{equation*}
we then show that ($dd^c(e^t)$ is bounded and that) $(d^cf-fd^ct)\wedge \omega^{m-1}$ restricted to the slice $\{t= s\}$ is a $O(e^{-s})$, with respect to $\eta\wedge p^*(\omega_{D_j})^{m-1}$ say.  \cqfd

~

\section{Generalization: the simple normal crossings case}  \label{part_gen}

Assume now $D$ has simple normal crossings. 
What follows is devoted to provide the necessary adaptation to get Theorem \ref{thm_topobstr} in this general case. 

We want first to construct circle actions around each component in order to have decompositions similar to \eqref{eqdecomp1}. 
For this we need the action around a fixed component to respect the other components. 
If these actions come from the identification of the normal bundles to tubular neighbourhoods with respect to some exponential map, we need the connection giving this exponential map to provide some orthogonality to intersecting component. 
In other words, we are asking for a connection with respect to which the crossing components are totally geodesic. 
We are also asking for the connection to preserve the complex structure $J$ of $X$, because we also need asymptotics on such an operator. 

There are obstructions to asking for such a connection to be the Levi-Civita connection of a smooth Kähler metric defined on some neighbourhood of the divisor. 
Nevertheless we only need the connection, not the metric, and we can construct it as follows. 

Take $\mathcal{U}=\mathcal{U}_1\sqcup\cdots\sqcup\mathcal{U}_m$ a finite open cover of $D$ in$X$, such that $U\in \mathcal{U}_k$ if, and only if, $U$ is an open set of holomorphic coordinates meeting a codimension $k$ crossing of $D$, and not meeting any codimension $\geq k+1$ crossing, this for $k=1,\dots,m$. 
Now if a codimension $k$ crossing, $D_1\cap\cdots\cap D_k$ say, is given by $\{z_1=\cdots= z_k=0\}$ in $U\in\mathcal{U}_k$, set $\nabla_U$ for the pull-back of the Levi-Civita of the euclidian metric in $U$. 
Take moreover a partition of unity $(\chi_U)_{U\in\mathcal{U}}$ associated to $\mathcal{U}$, and set 
 \begin{equation*} \label{eq_df_nabla}
  \nabla =\sum_{U\in\mathcal{U}}\chi_U\nabla_U
 \end{equation*}
on the neighbourhood $\bigcup_{U\in\mathcal{U}}U$ of $D$. 

It is clear that
 \begin{itemize}
  \item[$\bullet$] each intersection made from the $D_j$ is totally geodesic for $\nabla$;
  \item[$\bullet$] the complex structure $J$ of $X$ is parallel for $\nabla$.
 \end{itemize}

Fix now $j\in\{1,\dots,N\}$, and set $D'_j=D_j\backslash\bigcup_{j'\neq j}D_{j'}$ (resp. $D'=\sum_{j\neq j'}D_{j'}$). 
Thanks to the exponential map of $\nabla$, one identifies a neighbourhood $\mathcal{V}$ of $D_j$ in $X$ (resp. neighbourhood $\mathcal{V}'$ of $D_j'$ in $X\backslash D'$) to the normal bundle of $D_j$ in $X$ (resp. to its restriction to $D_j'$). 
We call $p_j$ the associated projection on $\mathcal{V}$; one has: $p_j(\mathcal{V}')=D_j'$

Through this identification, one gets an $S^1$ action around $D_j$ preserving the other components of $D$, and the discs in the resulting fibration of $\mathcal{V}$ are asymptotically holomorphic. 

Averaging $|\sigma_j|_j^2$ as in part \ref{part_fibr}, one gets an $S^1$ invariant function $t_j$ verifying $t_j=u_j+o(1)$ at any order with respect to a Poincaré metric on $X\backslash D_j$, $\omega_j$ say (e.g. $\omega_j=\omega_0-dd^c u_j$). 
One also builds a connection 1-form $\eta_j$ giving length $2\pi$ to the fibers of the action. 
Finally, if $g$ is the Riemannian metric associated to $\omega$ and $h_j$ that of $\omega|_{D_j'}$ (the Kähler form of Poincaré type induced on $D_j'$), one has on $\mathcal{V}'$ the asymptotics 
 \begin{equation*} \label{eq_asympg_snc}
  g=dt_j^2 + 4e^{-t_j}\eta_j^2+p_j^*h_j+o(1)
 \end{equation*}
with the perturbation $o(1)$ understood at any order with respect to $\omega_j$. 

Similarly, when it makes sense (and it does for Poincaré type potentials, for example), define $\Pi_{0}f$ as the $S^1$-invariant part of $f$ around $D'_j$, and $f_{0,j}$ and $f_{1,j}$ by 
 \begin{equation*}
  f_{0,j}(t_j)=\frac{1}{\vl(D_j)}\int_{D'_j}(p_j)_*(\Pi_{0}f)(t_j,\cdot)\vol^{\omega|_{D'_j}} \quad 
   \text{and} \quad 
  f_{1,j}=\Pi_{0}f-f_{0,j}
 \end{equation*}
(here $\vl(D_j)$ is computed with respect to $\omega|_{D_j'}$, and is thus equal to $\tfrac{[\omega_0|_{D_j}]^{m-1}}{(m-1)!}$). 
Derivation formulas \eqref{eq_deriv2} and \eqref{eq_deriv3} become
 \begin{equation*} \label{eq_deriv2_snc}
  d^cf = Jdf = 2\big(\partial_{t_j} f_0(t)+\partial_{t_j} f_1(t,z)\big)e^{-t_j}\eta_j+d_{D'_j}^cf_1(t,z)+o(1)
 \end{equation*}
with $d_{D'_j}^c=p_j^*(J_{D_j}d)$, and
 \begin{equation*}  \label{eq_deriv3_snc}
  \begin{aligned}
   dd^cf  = & 2\big(\partial_{t_j}^2 f_0(t_j)+\partial_{t_j}^2 f_1(t_j,z)
                    -\partial_{t_j} f_0(t_j)-\partial_{t_j} f_1(t_j,z)\big)e^{-t_j}dt_j\wedge\eta_j \\
            & + 2e^{-t_j}d_{D'_j}\big(\partial_{t_j} f_0(t_j)+\partial_{t_j} f_1(t_j,z)\big)\wedge\eta_j 
              +dt_j\wedge d_{D'_j}^c\partial_{t_j}f_1(t_j,z) \\
            & + dd^c_{D'_j}f_1(t_j,z)+o(1)
  \end{aligned}
 \end{equation*}
for $f\in C^{k,\alpha}(\XD)$, $k\geq2, \alpha\in(0,1)$, or $f\in \tmom$ (in both cases, the derivatives of $f$ with respect to $\omega$ and of order $\leq k-1$ are bounded, and the $S^1$ depending part is included in the $o(1)$). 

The point to be noticed here is that these estimations are uniform along $D'_j$, away from the other $D_{j'}$ as well as close to them.

Now, a careful reading of the proofs of Propositions \ref{prop_calcul1} and \ref{prop_calcul2} allows us to state them in our normal crossings case, with small adaptations: if $\varphi\in\tmom$, then 
 \begin{equation} \label{eq_calcul1_snc}
  \int_{\{t_j\leq s\}}e^{t_j}\omega_{\varphi}^m = 4\pi m!\vl(D_j)\big(s-\varphi_{0,j}(s)\big)+O(1), 
 \end{equation}
with moreover this integral tending to $+\infty$ when $s$ goes to $+\infty$, and 
 \begin{equation}  \label{eq_calcul2_snc}
  \int_{\{t_j\leq s\}}e^{t_j}\scal_{\varphi}\omega_{\varphi}^m = 4\pi m!
                                                        \vl(D_j)\big((\overline{\scal}_{D_j}-2)s-\overline{\scal_{D'_j}}\varphi_{0,j}(s)\big)+O(1)
 \end{equation}
where now $\overline{\scal}_{D_j}$ is computed for the class of \textit{Poincaré type} Kähler metrics on $D'_j$ relatively to $\omega_0|_{D_j}$ (instead of smooth metrics in the smooth divisor case), and hence is given by
 \begin{equation*}
  \overline{\scal_{D_j'}} = -4\pi(m-1)\frac{c_1\big(K[D]|_{D_j}\big)[\omega_0|_{D_j}]^{m-2}}{[\omega_0|_{D_j}]^{m-1}}
                          = -4\pi(m-1)\frac{c_1\big([D_j]\big)c_1\big(K[D]\big)[\omega_0]^{m-1}}{c_1\big([D_j]\big)[\omega_0]^{m-1}}.
 \end{equation*}
For example, Stokes theorem, which is a key-ingredient in those proofs, work again since the different integrands (with fixed $s$) in the integration by parts are $L^1$ with respect to $\omega|_{\{t_j\leq s\}}$. 
Moreover the $O(e^{-s})$ of those proofs were much better that needed, and the present $o(1)$ at our disposal are sufficient to get the final $O(1)$.

Assuming finally that $\omega_{\varphi}$ has constant scalar curvature and comparing $\overline{\scal}$ times \eqref{eq_calcul1_snc} with \eqref{eq_calcul2_snc} 
gives us that $\overline{\scal}_{D_j}>\overline{\scal}$. 
This holds for all $j=1,\dots,N$, and thus Theorem \ref{thm_topobstr} is proved in the simple normal crossings case.
~

\begin{footnotesize}

 \renewcommand{\refname}{References}

\end{footnotesize}

~

\small \textsc{École Normale Supérieure, UMR 8553}

\end{document}